\documentclass[11pt]{amsart}
\usepackage{amsmath,amsthm,amsfonts,amssymb,latexsym,epsfig}

\newtheorem{theorem}{Theorem}[section]

\numberwithin{equation}{section}

\theoremstyle{definition}

\theoremstyle{remark}

\begin{document}
\title{A note on $l^2$ norms of weighted mean matrices}
\author{Peng Gao}
\address{Department of Computer and Mathematical Sciences,
University of Toronto at Scarborough, 1265 Military Trail, Toronto
Ontario, Canada M1C 1A4} \email{penggao@utsc.utoronto.ca}
\date{June 11, 2007.}


\begin{abstract}
 We give a proof of Cartlidge's result on the $l^{p}$ operator norms of weighted mean
matrices for $p=2$ on interpreting the norms as eigenvalues of
certain matrices.
\end{abstract}

\maketitle
\section{Introduction}
\label{sec 1} \setcounter{equation}{0}

  Suppose throughout that $p\neq 0, \frac{1}{p}+\frac{1}{q}=1$.
   Let $l^p$ be the Banach space of all complex sequences ${\bf a}=(a_n)_{n \geq 1}$ with norm
\begin{equation*}
   ||{\bf a}||: =(\sum_{n=1}^{\infty}|a_n|^p)^{1/p} < \infty.
\end{equation*}
  The celebrated
   Hardy's inequality (\cite[Theorem 326]{HLP}) asserts that for $p>1$,
\begin{equation}
\label{eq:1} \sum^{\infty}_{n=1}\Big{|}\frac {1}{n}
\sum^n_{k=1}a_k\Big{|}^p \leq (\frac
{p}{p-1})^p\sum^\infty_{k=1}|a_k|^p.
\end{equation}
   Hardy's inequality can be regarded as a special case of the
   following inequality:
\begin{equation*}
\label{01}
   \sum^{\infty}_{j=1}\big{|}\sum^{\infty}_{k=1}c_{j,k}a_k
   \big{|}^p \leq U \sum^{\infty}_{k=1}|a_k|^p,
\end{equation*}
   in which $C=(c_{j,k})$ and the parameter $p$ are assumed
   fixed ($p>1$), and the estimate is to hold for all complex
   sequences ${\bf a}$. The $l^{p}$ operator norm of $C$ is
   then defined as the $p$-th root of the smallest value of the
   constant $U$:
\begin{equation*}
\label{02}
    ||C||_{p,p}=U^{\frac {1}{p}}.
\end{equation*}

    Hardy's inequality thus asserts that the Ces\'aro matrix
    operator $C$, given by $c_{j,k}=1/j , k\leq j$ and $0$
    otherwise, is bounded on {\it $l^p$} and has norm $\leq
    p/(p-1)$. (The norm is in fact $p/(p-1)$.)

    We say a matrix $A$ is a summability matrix if its entries satisfy:
    $a_{j,k} \geq 0$, $a_{j,k}=0$ for $k>j$ and
    $\sum^j_{k=1}a_{j,k}=1$. We say a summability matrix $A$ is a weighted
    mean matrix if its entries satisfy:
\begin{equation}
\label{021}
    a_{j,k}=\lambda_k/\Lambda_j,  ~~ 1 \leq k \leq
    j; \Lambda_j=\sum^j_{i=1}\lambda_i, \lambda_i \geq 0, \lambda_1>0.
\end{equation}

    Hardy's inequality \eqref{eq:1} now motivates one to
    determine the $l^{p}$ operator norm of an arbitrary summability matrix $A$.
   In an unpublished dissertation \cite{Car}, Cartlidge studied
weighted mean matrices as operators on $l^p$ and obtained the
following result (see also \cite[p. 416, Theorem C]{B1}).
\begin{theorem}
\label{thm02}
    Let $1<p<\infty$ be fixed. Let $A$ be a weighted mean matrix given by
    \eqref{021}. If
\begin{equation}
\label{022}
    L=\sup_n(\frac {\Lambda_{n+1}}{\lambda_{n+1}}-\frac
    {\Lambda_n}{\lambda_n}) < p ~~,
\end{equation}
    then
    $||A||_{p,p} \leq p/(p-L)$.
\end{theorem}

   We note here there are several published proofs of Cartlidge's
   result. Borwein \cite{Bor} proved a far more general
result than Theorem \ref{thm02} on the $l^p$ norms of generalized
Hausdorff matrices. Rhoades \cite[Theorem 1]{Rh} obtained a slightly
general result
 than Theorem \ref{thm02}, using a modification of the proof of
 Cartlidge. Recently, the author \cite{G1} also gave a simple proof of Theorem
 \ref{thm02}.

   It is our goal in this note to give another proof of Theorem \ref{thm02} for the case $p=2$,
   following an approach of Wang and Yuan in \cite{W&Y}, which
   interprets the left-hand side of \eqref{eq:1} when $p=2$ as a quadratic form
   so that Hardy's inequality follows from estimations of the
   eigenvalues of the corresponding matrix associated to the quadratic form.
   We will show in the next section that the same idea also works for the case of
   weighted mean matrices.

\section{Proof of Theorem \ref{thm02} for $p=2$}
\label{sec 2} \setcounter{equation}{0}
   We may assume $a_n$ being real without loss of generality and
   it suffices to prove the theorem for any finite summation from $n=1$ to $N$ with $N \geq 1$.
   We also note that it follows from our assumption on $L$ that $\lambda_n>0$. Now consider
\begin{equation*}
   \sum^{N}_{n=1}\Big ( \sum^{n}_{i=1}\frac {\lambda_i}{\Lambda_n}a_i\Big )^2
   =\sum^{N}_{n=1}\Big ( \sum^{n}_{i,j=1}\frac {\lambda_i\lambda_j}{\Lambda^2_n}a_ia_j\Big )
   =\sum^{N}_{n=1}\alpha_{i,j}a_ia_j, \hspace{0.1in}
   \alpha_{i,j}= \sum^{N}_{k \geq \max{(i,j)}}\frac {\lambda_i\lambda_j}{\Lambda^2_k}.
\end{equation*}
  We view the above as a quadratic form and define the associated
  matrix $A$ to be
\begin{equation*}
  A=\Big ( \alpha_{i, j} \Big )_{1 \leq i, j \leq N}.
\end{equation*}
   We note that the matrix $A$ here is certainly positive definite, being equal to $B^{t}B$
   with $B$ a lower-triangular matrix,
\begin{equation*}
   B=\Big ( b_{i, j} \Big )_{1 \leq i, j \leq N},  \hspace{0.1in} b_{i,j}=\lambda_j/\Lambda_i,  ~~ 1 \leq j \leq
    i;   \hspace{0.1in} b_{i,j}=0, j > i.
\end{equation*}
   It is easy to check that the entries of $B^{-1}$ are given by
\begin{equation*}
  \big (  B^{-1} \big )_{i,i}=\frac
  {\Lambda_i}{\lambda_i}, \hspace{0.1in}
  \big (  B^{-1} \big )_{i+1,i}=-\frac
  {\Lambda_i}{\lambda_{i+1}}, \hspace{0.1in} \big (  B^{-1} \big
  )_{i,j}=0 \hspace{0.1in} \text{otherwise}.
\end{equation*}

    In order to
   establish our assertion, it suffices to show that the maximum eigenvalue of $A$ is less than
   $4/(2-L)^2$ or the minimum eigenvalue of its inverse $A^{-1}$ is greater than
   $(2-L)^2/4$ which is equivalent to proving that the matrix $A^{-1} - \lambda I_N$ is positive
   definite, where $\lambda=(2-L)^2/4$ and $I_N$ is the $N \times N$ identity matrix. Using the expression $A^{-1}=B^{-1}(B^{-1})^{t}$, we
   see that this is equivalent to showing that 
   for any integer $N \geq 1$ and any
   real sequence ${\bf a}=(a_n)_{1 \leq n \leq N}$,
\begin{equation}
\label{3.0}
   \sum^{N-1}_{n=1}\Big ( \frac
  {\Lambda_n}{\lambda_n}a_n-\frac
  {\Lambda_n}{\lambda_{n+1}}a_{n+1}\Big )^2 +  \frac
  {\Lambda^2_N}{\lambda^2_N}a^2_N\geq \frac {(2-L)^2}{4}\sum^N_{n=1}a_n^2.
\end{equation}
  
  For any integer $n \geq 1$ and fixed constants $\alpha$, $\beta$, $a_{n+1}, \mu_n$ (here $\alpha$, $\beta$ may depend on $n$), we consider the following function:
\begin{equation*}
  f(a_n) :=(\alpha a_n-\beta a_{n+1})^2-\mu_na^2_n.
\end{equation*}
   When $\mu_n > \alpha^2$, it is easy to see that 
\begin{equation}
\label{2.2}
  f(a_n) \leq f(\frac {\alpha \beta a_{n+1}}{\alpha^2-\mu_n})=\frac {\beta^2\mu_{n}\alpha^2_{n+1}}{\mu_n-\alpha^2},
\end{equation}
  with the above inequality reversed when $\mu_n < \alpha^2$.

  On taking $\alpha=\Lambda_n/\lambda_n, \beta=\Lambda_n/\lambda_{n+1}$ here, we obtain that for any $0 < \mu_n < \Lambda^2_n/\lambda^2_n$,
\begin{equation*}
   \Big ( \frac
  {\Lambda_n}{\lambda_n}a_n-\frac
  {\Lambda_n}{\lambda_{n+1}}a_{n+1}\Big )^2- \mu_n a^2_n \geq -\frac {\Lambda^2_n/\lambda^2_{n+1}}
  {\Lambda^2_n/\lambda^2_n-\mu_n}\mu_na^2_{n+1}.
\end{equation*}
   Summing the above inequality for $n=1, \ldots, N-1$ yields:
\begin{eqnarray}
\label{3.3}
 &&  \sum^{N-1}_{n=1}\Big ( \frac
  {\Lambda_n}{\lambda_n}a_n-\frac
  {\Lambda_n}{\lambda_{n+1}}a_{n+1}\Big )^2 +  \frac
  {\Lambda^2_N}{\lambda^2_N}a^2_N \\
  &\geq &  \mu_1a^2_1+\sum^{N-2}_{n=1}\Big(\mu_{n+1}- \frac {\Lambda^2_n/\lambda^2_{n+1}}
  {\Lambda^2_n/\lambda^2_n-\mu_n}\mu_n\Big )a_{n+1}^2+ \Big ( \frac
  {\Lambda^2_N}{\lambda^2_N}  - \frac {\Lambda^2_{N-1}/\lambda^2_N}
  {\Lambda^2_{N-1}/\lambda^2_{N-1}-\mu_{N-1}}\mu_{N-1} \Big ) a^2_N.
  \nonumber
\end{eqnarray}
   We now want to find a number $\mu _n$ satisfying  $0 < \mu_n < \Lambda^2_n/\lambda^2_n$ for any integer $n \geq 1$, 
   such that the following inequality holds for $1 \leq n \leq N-1$:
\begin{equation}
\label{3.4}
  \mu_{n+1}- \frac {\Lambda^2_n/\lambda^2_{n+1}}
  {\Lambda^2_n/\lambda^2_n-\mu_n}\mu_n \geq \frac {(2-L)^2}{4}.
\end{equation}

   For this purpose, we set
\begin{equation*}
  \mu_n=(k+c)\frac {\Lambda_n}{\lambda_n}-c, \hspace{0.1in} k=(2-L)^2/4,
\end{equation*}
     with $c$ a constant to be specified later.  For the so chosen $\mu_n$'s, inequality \eqref{3.4} can be seen to be equivalent to
\begin{equation}
\label{3.6'}
    (k+c)x^2+(k+c)x-\big ( (k+c)x-c \big )y-(k+c) \big ( (k+c)x-c \big ) -c \geq 0.
\end{equation}
   where $x=\Lambda_n/\lambda_n$ and $y= \Lambda_{n+1}/\lambda_{n+1}$. Note that by the assumption \eqref{022}, it follows that $y \leq x + L$ and the case $n=1$
of \eqref{022} implies $L>0$. Note also that $x \geq 1$ so that the left-hand side expression of inequality \eqref{3.6'} is an decreasing function of $y$ for fixed $x$. Hence we can replace
   $y$ there and conclude that \eqref{3.6'} follows from the following inequality:
\begin{equation*}
    (k+c)x^2+(k+c)x-\big ( (k+c)x-c \big )(L+x)-(k+c) \big ( (k+c)x-c \big ) -c \geq 0.
\end{equation*}
    Equivalently, we can recast the above inequality as:
\begin{equation}
\label{2.5}
    \Big(k+2c-L(k+c)-(k+c)^2 \Big)x+cL+c(k+c)-c \geq 0.
\end{equation}  
   In order for the above inequality to hold for all $x \geq 1$, we need to choose $c$ so that
\begin{equation}
\label{2.6}
    k+2c-L(k+c)-(k+c)^2 \geq 0.
\end{equation}  
   We now choose the value of $c$ so that the left-hand side expression above when considered as a function of $c$ is maximized. It is easy to see that in this case, $c+k=(2-L)/2$ so that
\begin{equation*}
  c=\frac {L(2-L)}{4}.
\end{equation*}
   It is then easy to check that for the so chosen $c$, the left-hand side expression of inequality \eqref{2.6} is reduced to $0$ and the left-hand side expression of inequality \eqref{2.5} becomes $cL/2 >0$. It follows that inequality \eqref{3.4} holds for $1 \leq n \leq N-1$. Note also that in our case $\mu_1=k$ and for $1 \leq n \leq N$,
\begin{equation*}
  0< \mu_n =(k+c)\frac {\Lambda_n}{\lambda_n}-c < (k+c)\frac {\Lambda_n}{\lambda_n} = \Big ( \frac {2-L}{2} \Big ) \frac {\Lambda_n}{\lambda_n} < \frac {\Lambda_n}{\lambda_n} \leq \frac
  {\Lambda^2_n}{\lambda^2_n},
\end{equation*}
   so that our assumption $0 < \mu_n < \Lambda^2_n/\lambda^2_n$ is satisfied. In particular, we have 
\begin{equation*}
   \frac
  {\Lambda^2_N}{\lambda^2_N}  - \frac {\Lambda^2_{N-1}/\lambda^2_N}
  {\Lambda^2_{N-1}/\lambda^2_{N-1}-\mu_{N-1}}\mu_{N-1} \geq  \mu_N- \frac {\Lambda^2_{N-1}/\lambda^2_N}
  {\Lambda^2_{N-1}/\lambda^2_{N-1}-\mu_{N-1}}\mu_{N-1} \geq \frac {(2-L)^2}{4}.
\end{equation*}
  From this we see that inequality \eqref{3.0} follows from inequality \eqref{3.3} and this completes the proof.
\section{Further Discussions}
\label{sec 3} \setcounter{equation}{0}
   We poin out here inequality \eqref{3.0} can be regarded as an analogue to the following discrete inequality of Wirtinger's type studied by Fan, Taussky and Todd \cite[Theorem 8]{FTT}:
\begin{equation}
\label{3.1}
  a^2_1+\sum^{N-1}_{n=1}(a_n-a_{n+1})^2 + a^2_N \geq 2 \Big (1-\cos \frac {\pi}{N+1} \Big )\sum^{N}_{n=1} a^2_n.
\end{equation}
    Converses of the above inequality was found by Milovanovi\'c and
    Milovanovi\'c \cite{MM}:
\begin{equation}
\label{3.2}
  a^2_1+\sum^{N-1}_{n=1}(a_n-a_{n+1})^2 + a^2_N \leq 2 \Big (1+\cos \frac {\pi}{N+1} \Big )\sum^{N}_{n=1} a^2_n.
\end{equation}
   Simple proofs of inequalities \eqref{3.1} and \eqref{3.2} were given by Redheffer \cite{R2} and Alzer \cite{A1}, respectively. Our proof of Theorem \ref{thm02} for $p=2$ in the previous section is motivated by the methods used in  \cite{R2} and \cite{A1}.
    
  To end this paper, we note the paper \cite{Lo} contains several generalizations of inequalities of \eqref{3.1} and \eqref{3.2}, one of them can be stated as:
\begin{theorem}
\label{thm3.1}
  For any real sequence ${\bf a}=(a_n)_{1 \leq n \leq N}$, and two positive real numbers $a, b$, 
\begin{equation}
\label{3.7}
  (a^2+b^2-2ab\cos \frac {\pi}{N+1} \Big )\sum^{N}_{n=1} a^2_n \leq b^2a^2_1+\sum^{N-1}_{n=1}(aa_n-ba_{n+1})^2 + a^2a^2_N \leq  \Big (a^2+b^2+2ab\cos \frac {\pi}{N+1} \Big )\sum^{N}_{n=1} a^2_n.
\end{equation}
\end{theorem}

  The proof given in \cite{Lo} to the above theorem is to regard
\begin{equation*}
  b^2a^2_1+\sum^{N-1}_{n=1}(aa_n-ba_{n+1})^2 + a^2a^2_N
\end{equation*}
   as a quadratic form with the associated
  matrix $A$ being symmetric tridiagonal with its entries
   given by
\begin{equation*}
  \big (  A \big )_{i,i}=a^2+b^2, \hspace{0.1in}
  \big (  A \big )_{i,i+1}=\big ( A \big )_{i+1,i}=-ab, \hspace{0.1in} \big (  A \big
  )_{i,j}=0 \hspace{0.1in} \text{otherwise}.
\end{equation*}
   The eigenvalues of $A$ are shown in \cite{Lo}  to be $a^2+b^2+2ab\cos(\frac {k\pi}{N+1}), 1 \leq k \leq N$, from which Theorem \ref{thm3.1} follows easily.

   We note here one can also give a proof of Theorem \ref{thm3.1} following the methods in  \cite{R2} and \cite{A1} as one checks readily that the right-hand side inequality of \eqref{3.7} follows on taking $\alpha=a, \beta=b, \mu_n=a^2+ab\sin(n+1)t/\sin(nt), t = \pi/(N+1)$ in inequality \eqref{2.2} and summing for $n=1, \ldots, N-1$.  Similarly, 
  the left-hand side inequality of \eqref{3.7} follows from on taking $\alpha=a, \beta=b, \mu_n=a^2-ab\sin(n+1)t/\sin(nt), t = \pi/(N+1)$ in inequality \eqref{2.2} (with inequality reversed there).

\end{document}